\documentclass[a4paper, 12pt]{scrartcl}

\usepackage[utf8]{inputenc}		
\usepackage{amsmath}				
\usepackage{mathtools}			
\usepackage{amsthm,amssymb,wasysym}		
\usepackage{amsfonts}			
\usepackage[numbers,square]{natbib}	
\usepackage{comment}			
\usepackage{authblk}			
\usepackage[babel]{csquotes}

\newcommand{\xbar}{\bar{x}}
\newcommand{\xring}{\mathring{x}}
\newcommand{\ybar}{\bar{y}}
\newcommand{\dxbar}{\overline{\dx}}
\newcommand{\dybar}{\overline{\dy}}
\newcommand{\dzbar}{\overline{\dz}}
\newcommand{\drbar}{\overline{\dr}}
\newcommand{\dxbarbar}{\overline{\overline{\dx}}}
\newcommand{\dVbar}{\overline{\dV}}
\newcommand{\dxi}{\dx_i}
\newcommand{\dyi}{\dy_i}
\newcommand{\dxbari}{\dxbar_i}
\newcommand{\dybari}{\dybar_i}
\newcommand{\dxbarbari}{\dxbarbar_i}
\newcommand{\Deltae}{\Delta e}
\newcommand{\Deltaei}{\Delta e_i}
\newcommand{\pdd}[2]{\frac{\partial #1}{\partial #2}}

\newcommand{\id}{ {1\!\!\!\:1 } }

\DeclareMathOperator{\dev}{dev}

\DeclareMathOperator{\tr}{tr}

\newcommand{\intd}[1]{\mathrm{d#1}}
\newcommand{\dr}{\mathrm{dr}}

\newcommand{\dx}{\mathrm{dx}}
\newcommand{\dy}{\mathrm{dy}}
\newcommand{\dz}{\mathrm{dz}}
\newcommand{\dV}{\mathrm{dV}}
\newcommand{\dphi}{\mathrm{d\varphi}}

\newcommand{\rbar}{\bar{r}}

\newcommand{\zbar}{\bar{z}}

\newcommand{\Hcomment}[1]{\footnote{#1}}
\newcommand{\Hgerman}[1]{\footnote{German: #1}}
\newcommand{\Herror}[1]{\footnote{A typographical error was corrected: #1}}
\newcommand{\Href}[1]{\footref{#1}}

\newcommand{\afrac}[2]{#1\! /\! #2}
\newcommand{\tel}[1]{\frac{1}{#1}}
\newcommand{\half}{\tel{2}}
\newcommand{\eps}{\varepsilon}
\newcommand{\norm}[1]{\Vert #1 \Vert}

\newcommand{\beginsmallpmatrix}{\left( \begin{smallmatrix}}

\newcommand{\nnl}{\nonumber \\}
\newcommand{\inv}{^{-1}}

\usepackage{titlesec}
\titleformat{\section}%
  [hang]
  {\normalfont\bfseries\Large}
  {}
  {0pt}
  {}

\setkomafont{dictumtext}{\normalfont\normalcolor\itshape\setlength{\parindent}{1em}\noindent}

\begin{document}

\title{The axiomatic deduction of the quadratic Hencky strain energy by Heinrich Hencky}
\author{Heinrich Hencky\\\vspace*{4mm}{\small translated and annotated by}\vspace*{2mm}\\Patrizio Neff\thanks{Corresponding author.\: Head of Chair for Nonlinear Analysis and Modelling, University of Duisburg-Essen, Thea-Leymann-Str. 9, 45127 Essen, Germany, email: patrizio.neff@uni-due.de},
Bernhard Eidel\thanks{Institute of Mechanics, University of Duisburg-Essen, Universitätsstr. 15, 45141 Essen, Germany, email: bernhard.eidel@uni-due.de}
and Robert Martin\thanks{Chair for Nonlinear Analysis and Modelling, University of Duisburg-Essen, Thea-Leymann-Str. 9, 45127 Essen, Germany, email: robert.martin@stud.uni-due.de }}
\date{\today}

\maketitle

\abstract{
The introduction of the quadratic Hencky strain energy based on the logarithmic strain tensor $\log V$ is a milestone in the development of nonlinear elasticity theory in the first half of the 20\textsuperscript{th} century. Since the original manuscripts are written in German, they are not easily accessible today. However, we believe that the deductive approach taken by Hencky deserves to be rediscovered today.\\
In this work we have gathered parts of the original contributions \enquote{Über die Form des Elastizitätsgesetzes bei ideal elastischen Stoffen}, \enquote{Welche Umstände bedingen die Verfestigung bei der bildsamen Verformung von festen isotropen Körpern?} and \enquote{Das Superpositionsgesetz eines endlich deformierten relaxationsfähigen elastischen Kontinuums und seine Bedeutung für eine exakte Ableitung der Gleichungen für die zähe Flüssigkeit in der Eulerschen Form} which center around this deductive approach. We tried to provide, for the first time, a faithful translation into English.\\
All footnotes are our addition.
}
\newpage

\section[Heinrich Hencky: On the form of the law of elasticity for ideally elastic materials]{Heinrich Hencky: On the form of the law of elasticity for ideally elastic materials\Hgerman{Über die Form des Elastizitätsgesetzes bei ideal elastischen Stoffen \cite{Hencky1928}\phantom{\cite{Hencky1929}\phantom{\cite{Hencky1929super}}}}}

\dictum{The law of elasticity for finite deformations. How do different states of stress and deformation superimpose in simple cases? The expression for elastic work and the theory of tension and compression of rubberlike elastic bodies.}
\subsection*{Introduction}
We call \emph{ideally elastic} a material which does not lose stored energy under arbitrarily large deformations and thus returns to its original state after unloading. It is not important that such an idealized elastic cycle does not actually exist and our ideally elastic material must therefore remain an ideal. Like so many mathematical and geometric concepts, it is a useful ideal, because once its deducible properties are known it can be used as a comparative rule for assessing the actual elastic behaviour of physical bodies.\\
The theory of elastic deformations, which has already been developed to a certain degree, will be useful for this task \cite[p.563-576]{hamel1965elementare}. Only one important aspect of the classical theory of finite deformations (c.f. \cite[p. 51-54]{1901encyklopaedie} as well as the references there) needs to be reconsidered: although it is correct to assume that the elastic energy is given through a function depending only on the rotational invariants\Hcomment{The rotational invariants are the principal invariants of the Finger tensor $B=FF^T$ or the right Cauchy-Green tensor $C=F^TF$.} of the deformation tensor, it is not appropriate to base the entire approach to the problem on this formal insight.\\
We will first show that it is possible to provide the law of elasticity in a direct way free of any arbitrariness. Furthermore we will see that, even in simple cases, the elastic energy is governed by laws too complex to be approached directly, in contrast to infinitesimal deformations.

\subsection{The law of elasticity for the ideally elastic body}
If we consider a volume element whose position and orientation are determined by the deformation, we can interpret its change under a deformation as a \emph{dilation} of the three mutually orthogonal edges of the volume element, a \emph{translation} of the center of mass and a finite \emph{rotation} about a particular \emph{axis}. This holds for finite deformations as well. However, there is an essential difference to the infinitesimal case, where the particle must actually pass through the specified changes of state. In finite deformations, only the final state is considered, while the intermediary states are not. This difference is especially important to the computation of the elastic energy.\\
While the edges of the parallelepiped mentioned above are pairwise orthogonal before and after the deformation, they are rotated against their original orientation by a finite angle. We denote these three principal directions of strain with the indices $1,~2$ and $3$ here and throughout. Note that those directions refer to a material particle and not a fixed spatial direction.\\
Furthermore, we denote by $\dx$ the length of the arc element in the undeformed state and by $\dxbar$ the length of the same element after the deformation, while $\dV$ and $\dVbar$ denote the corresponding volume elements. We find
\[
	\dV = \dx_1\,\dx_2\,\dx_3 \quad \text{ and } \quad \dVbar = \dxbar_1\,\dxbar_2\,\dxbar_3\,.
\]
The three principal strains are defined by the following equalities:\Hcomment{Hencky introduces what is known today as Swainger's strain tensor $\id-V\inv$, where $V=\sqrt{FF^T}$ is the left Biot stretch tensor.}
\begin{gather*}
	e_1 = \frac{\dxbar_1 - \dx_1}{\dxbar_1}\,; \qquad e_2 = \frac{\dxbar_2 - \dx_2}{\dxbar_2}\,; \\
	e_3 = \frac{\dxbar_3 - \dx_3}{\dxbar_3}\,.
\end{gather*}
To shorten notation we will use the index $i$ for all three principal directions, meaning $i$ is one of the indices $1,2,3$. Then we can write the above definitions of the strain as
\begin{subequations}
\begin{equation}
	e_i = \frac{\dxbar_i - \dx_i}{\dxbar_i}\,.
\end{equation}
We can now easily compute the ratio of the deformed length to the initial length\Hcomment{The quantities given by $\lambda_i=\frac{\dxbari}{\dxi} = \tel{1-e_i}$ are the \emph{principal stretches}, i.e. the eigenvalues of $V$. Therefore $e_i = 1-\frac1{\lambda_i}$ are the eigenvalues of $\id-V\inv$.},
\begin{equation}
	\frac{\dxbari}{\dxi} = \tel{1-e_i}\,,
\end{equation}
while the ratio of the volume in the deformed state to the initial volume is\Hcomment{Note that $\frac{\dVbar}{\dV} = \det F$ for the deformation gradient $F$.}
\begin{equation}
	\frac{\dVbar}{\dV} = \tel{(1-e_1)\,(1-e_2)\,(1-e_3)}\,.
\end{equation}
\end{subequations}
If we consider an element, cut out of the deformed body along the principal directions of strain in the final state of equilibrium, we must apply the three principal stresses\Hcomment{The stresses $S_i$ are the principal Cauchy stresses.} $S_1$, $S_2$ and $S_3$ to the principal directions. Of course those stresses refer to the final state. While in the field of material testing one might consider stresses referring to the initial state, we must, from the theoretical point of view, protest against such a conceptional abomination\Hgerman{Unbegriff}. We will give a relation between $S_i$ and $e_i$ such that the transition to the infinitesimal case yields Hooke's law in the form commonly used in engineering.\\
Let $G$ denote the shear modulus, $m$ the lateral contraction number\Hcomment{The \emph{Querkontraktionsziffer} (or \emph{Poisson number}) $m$ is the inverse of \emph{Poisson's ratio}: $m=\tel\nu$\,.}, which loses its meaning in the case of finite deformations, and define $k$ by\Hcomment{The parameter $k$ can also be written as $k=\frac{1+\nu}{3(1-2\nu)} = \frac{K}{2\,G}$, where $K$ is the \emph{bulk modulus}.}
\begin{subequations}
\begin{equation}
	k=\frac{m+1}{3\,(m-2)}\,.\label{eq:Hencky2a}
\end{equation}
By writing
\begin{equation}
	e = \tel3 \cdot(e_1+e_2+e_3)
\end{equation}
and
\begin{equation}
	S = \tel3 \cdot(S_1+S_2+S_3)
\end{equation}
we can state Hooke's law as
\begin{equation}
	S_i = 2\,G\,\{e_i+(3k-1)\cdot e\}\label{eq:Hencky2d}
\end{equation}
for the three principal stresses. Summation of these three equations yields
\begin{equation}
	S = 2\,G\,k\cdot 3\,e\,.\label{eq:Hencky2e}
\end{equation}
\end{subequations}
The left hand term of this equality is the hydrostatic part of the stress, which must be proportional to the change of volume in an isotropic material.\\
For finite stretches, however, the term $3e$ no longer denotes the change of volume, and thus the approach loses its mechanical meaning. Furthermore, while $e_i$ attaining the value $1$ means an infinite  deformation of the element, the corresponding stress does not become infinite, again contradicting mechanical reasoning. This contradiction can not be circumvented by choosing another definition of strain, as one can readily check.\\
If we want to preserve the simple form of the law of elasticity with only two elasticity parameters as well as to avoid the aforementioned contradictions, there is only one possible ansatz, which is the following:
\begin{subequations}
\label{eqs:Hencky3}
\begin{equation}
	S_i = 2\,G\cdot\ln\left\{ \frac{\dxbari}{\dxi}\cdot\left(\frac{\dVbar}{\dV}\right)^{k-\afrac13} \right\}\,.\label{eq:Hencky3a}
\end{equation}
By summation of the three principal stretches we obtain the hydrostatic part\Hcomment{The hydrostatic part is given by \[\tr S = 2\,G\,k \cdot\ln(\det F) = K\cdot \ln(\det F)\,,\] where $\tr S$ is the trace of the Cauchy stress tensor $S$.}, given by
\begin{equation}
	S = 2\,G\,k\cdot\ln\left\{\frac{\dVbar}{\dV}\right\}\,.\label{eq:Hencky3b}
\end{equation}
\end{subequations}
We can see that this equality relates the change of volume to the hydrostatic tension or compression in the correct way.\\
If the ratios inside the logarithms are close to $1$, the logarithms can be expanded into power series whose higher order terms can be omitted, transforming equations \eqref{eq:Hencky3a} and \eqref{eq:Hencky3b} into equations \eqref{eq:Hencky2d} and \eqref{eq:Hencky2e}.\\
For applications it is often more comfortable to introduce the strains $e_i$. We obtain
\begin{subequations}
\label{eqs:Hencky4}
\begin{align}
	S_i &= -2\,G\cdot\ln\left\{(1-e_i)\cdot((1-e_1)(1-e_2)(1-e_3))^{k-\afrac13}\right\}\,,\label{eq:Hencky4a}\\[2mm]
	S &= -2\,G\,k\ln\{(1-e_1)(1-e_2)(1-e_3)\}\,.\label{eq:Hencky4b}
\end{align}
\end{subequations}
The numbers $e_i$ can attain all values between negative infinity and $1$, including $0$. If $e_i$ is negative infinity or $1$ we obtain infinitely large values for the stress, just as required by the mechanical meaning of the law of elasticity. Thus equations \eqref{eqs:Hencky3} and \eqref{eqs:Hencky4} yield the desired law of elasticity for isotropic materials.

\subsection{Superimposing different states of strain and stress}
We will now consider the following question: what happens if we apply an additional load to an already stressed body, i.e. how do stresses and deformations superimpose? We can confine our considerations to the case of homogeneous stress and deformation, and we will distinguish three different states: State I is the non-loaded initial state, and if $\dx_1, \dx_2, \dx_3$ are lengths of the sides of a rectangular parallelepiped which remains rectangular after the first loading, ending in state II, this deformed parallelepiped will have lengths $\dxbar_1,\dxbar_2,\dxbar_3$. According to our definition of \;\emph{strain = final length - initial length : final length}\; we will find
\[
	e_i = \frac{\dxbar_i - \dx_i}{\dxbar_i}
\]
or, in a different notation,
\begin{subequations}
\begin{equation}
	\dxbar_i = \dx_i : (1-e_i)\,.
\end{equation}
If we now apply another deformation, our rectangular parallelepiped will not remain rectangular. However, we can find a different rectangular parallelepiped $\dy_1, \dy_2, \dy_3$ in state I which is oblique-angled in state II, but rectangular in state III with side lengths $\dybar_1, \dybar_2, \dybar_3$. The values of strain between states I and III are
\begin{equation}
	e_i' = \frac{\dybari - \dyi}{\dybari} \qquad \text{or} \qquad \dybari = \dyi : (1-e_i')\,.
\end{equation}
\end{subequations}
If we denote the principal stresses in state II by $S_i$ and those in state III by $S_i'$, the same relations hold between $S_i$ and $e_i$ on the one hand and $S_i'$ and $e_i'$ on the other hand, which are given by equations \eqref{eq:Hencky4a} and \eqref{eq:Hencky4b}. Thus if the lengths $\dxi$ and $\dyi$ are equal, we obtain
\begin{equation}
	S_i' - S_i = 2\,G\cdot\ln\left\{\frac{\dybari}{\dxbari}\cdot\left[\frac{\dVbar(\ybar)}{\dVbar(\xbar)}\right]^{k-\afrac13}\right\}\,.
\end{equation}
If $e_i$ and $e_i'$ are finite, the term inside the logarithm on the right hand side is \emph{not} the transformation of II to III, not even if the difference II$-$III is infinitesimal. Since the treatment of this general case would demand an effort disproportionate to the obtainable practical results, we will assume henceforth that the principal axes of stress and strain are parallel to the axes of a fixed coordinate system. Then the material parallelepiped remains constant and we obtain particularly simple rules for the superposition of different stresses and deformations as well as for the performed work.\\
The side lengths of the volume element are
\begin{align*}
	&\text{$\dx_1,\,\dx_2,\,\dx_3\:$ in state I,}\\
	&\text{$\dxbar_1,\,\dxbar_2,\,\dx_3\:$ in state II,}\\
	&\text{$\dxbarbar_1,\,\dxbarbar_2,\,\dxbarbar_3\:$ in state III,}
\end{align*}
while the transformations from \,I to II,\, I to III\, and \,II to III\, are
\begin{subequations}
\begin{equation}
	\left.\begin{alignedat}{2}
		&\dxbari = \dxi:(1-e_i) &&\:\text{ from I to II,}\\
		&\dxbarbari = \dxi:(1-e_i') &&\:\text{ from I to III,}\\
		&\dxbarbari = \dxbari:(1-\Delta e_i) &&\:\text{ from II to III.}
	\end{alignedat}\right\}\label{eq:Hencky7a}
\end{equation}
This implies
\begin{equation}
	1-\Deltaei = \frac{1-e_i'}{1-e_i}\,.\label{eq:Hencky7b}
\end{equation}
\end{subequations}
Subtracting $S_i$ from $S_i'$ yields
\[
	S_i' - S_i = \Delta S_i = -2\,G\cdot\ln\frac{(1-e_i')\,((1-e_1')(1-e_2')(1-e_3'))^{k-\afrac13}}{(1-e_i)\,((1-e_1)(1-e_2)(1-e_3))^{k-\afrac13}}
\]
and by using equality \eqref{eq:Hencky7b} we obtain
\begin{subequations}
\begin{equation}
	S_i' - S_i = \Delta S_i = -2\,G\cdot\ln[(1-\Deltaei)\cdot((1-\Deltae_1)(1-\Deltae_2)(1-\Deltae_3))^{1-\afrac13}]\,.\label{eqs:Hencky8a}
\end{equation}
Thus the difference $\Delta S_i$ depends only on the second transformation and is independent of the previous transformation. As mentioned above, this does not hold in the general case, even if the $\Delta S_i$ are infinitesimally small. In the latter case, the $\Deltaei$ are infinitesimal as well, and applying the series expansion of the logarithm to equation \eqref{eqs:Hencky8a} yields
\begin{equation}
	S_i' - S_i = \Delta S_i = 2\,G\cdot\{\Deltaei + (k-\afrac13)(\Deltae_1+\Deltae_2+\Deltae_3)\}\,,\label{eqs:Hencky8b}
\end{equation}
\end{subequations}
where the $\Delta S_i$ take the role of the stress tensors themselves. In this special case, and only in this special case, Hooke's law is the incremental law\Hcomment{Hooke's law is obtained as a first order approximation.} corresponding to our law of elasticity.

\subsection{The stored elastic energy}
The last sentence will be useful for computing the elastic work. Again we assume that the states II and III are infinitesimally close.\\
If $x_i$ are the coordinates of a single point in state I, its coordinates in states II and III are given by $x_i+u_i$ and $x_i+u_i+\Delta u_i$, respectively.\\
Then equations \eqref{eq:Hencky7a} take on the form
\begin{align*}
	\dxbarbari &= \dxi\left\{1+ \pdd{u_i}{x_i} + \pdd{\Delta u_i}{x_i}\right\}
\intertext{and}
	\dxbari &= \dxi\left\{1+\pdd{u_i}{x_i}\right\}\,,
\end{align*}
and comparison with equations \eqref{eq:Hencky7a} and \eqref{eq:Hencky7b} yields
\begin{align*}
	1-\Deltaei &= \frac{1+\pdd{u_i}{x_i}}{1+\pdd{u_i}{x_i}+\pdd{\Delta u_i}{x_i}}\,,\\
	1-\Deltaei &= \frac{1}{1+(1-e_i)\cdot\pdd{\Delta u_i}{x_i}}\,.
\end{align*}
From this we obtain the differential quotient of the additional displacement
\[
	\pdd{\Delta u_i}{x_i} = \frac{\Deltaei}{(1-\Deltaei)(1-e_i)}\,,
\]
and for an infinitesimal additional displacement we find
\begin{equation}
	\pdd{\Delta u_i}{x_i} = \frac{\Deltaei}{1-e_i}\,.\label{eqs:Hencky9}
\end{equation}
All occurring functions refer to the initial positions. We will denote the work with respect to the unit of the volume element in state I by $A_a$, while the work with respect to the volume element in state II will be denoted by $A_e$. Then the elastic energy of the element $\dx_1\,\dx_2\,\dx_3$ will increase by\Herror{The index $i$ was changed to $1$.}
\begin{align*}
	\Delta A_a \cdot \,\dx_1\,\dx_2\,\dx_3 &= S_1\cdot \,\dxbar_2\,\dxbar_3\,\cdot\,\pdd{\Delta u_1}{x_1}\,\cdot\,\dx_1\\
	&\quad + S_2\,\cdot\, \dxbar_3\,\dxbar_1\,\cdot\,\pdd{\Delta u_2}{x_2}\,\cdot\,\dx_2 + S_3\,\cdot\, \dxbar_1\,\dxbar_2\,\cdot\,\pdd{\Delta u_3}{x_3}\,\cdot\,\dx_3\,,
\end{align*}
and using equations \eqref{eq:Hencky7a} we obtain
\begin{align*}
	\Delta A_a = \tel{(1-e_1)(1-e_2)(1-e_3)}&\left\{ S_1(1-e_1)\cdot\pdd{\Delta u_1}{x_1}\right.\\
	&\quad \left.+ S_2(1-e_2)\cdot\pdd{\Delta u_2}{x_2} + S_3(1-e_3)\cdot\pdd{\Delta u_3}{x_3} \right\}\,.
\end{align*}
According to equation \eqref{eq:Hencky4b}, we find
\[
	(1-e_1)(1-e_2)(1-e_3) = e^{-\afrac{S}{K}}\,,
\]
where $K=2\,Gk$ denotes the elasticity modulus of uniform expansion or compression (on all sides)\Hcomment{The parameter $K$ denotes the \emph{bulk modulus}.}, and using equation \eqref{eqs:Hencky9} we can compute the increase of work:
\[
	\Delta A_a = e^{\afrac{S}{K}}\cdot\left[ S_1\Deltae_2+S_1\Deltae_2+S_3\Deltae_3 \right]\,.
\]
Equation \eqref{eqs:Hencky8b} allows us to change the last expression to an even more familiar form. We can invert these equations and, with\footnote{\emph{Young's modulus} is denoted by $E$.} $E=2\,G(1+\afrac1m)$, we obtain
\begin{align*}
	E\cdot\Deltae_1 &= \Delta S_1 - \tel{m}(\Delta S_2 + \Delta S_3)\,,\\
	E\cdot\Deltae_2 &= \Delta S_2 - \tel{m}(\Delta S_3 + \Delta S_1)\,,\\
	E\cdot\Deltae_3 &= \Delta S_3 - \tel{m}(\Delta S_1 + \Delta S_2)\,.\\
\end{align*}
We substitute the $\Delta$s by total differentials. It is easy to see that we obtain the familiar expression for work, since
\[
	S_1\,\intd{}e_1 + S_2\,\intd{}e_2 + S_3\,\intd{}e_3 = \half \cdot \intd{}J\,,
\]
where
\begin{subequations}
\begin{equation}
	J = S_1^2 + S_2^2 + S_3^2 - \frac2m(S_1S_2+S_3S_1+S_2S_3)\,.
\end{equation}
Thus we finally obtain the energy with respect to state I,
\begin{equation}
	A_a = \tel{2E}\cdot\int e^{\afrac{S}{K}}\cdot\intd{}J\,,\label{eq:Hencky10b}
\end{equation}
while the energy with respect to state II, i.e. to the unit volume in state II, computes to
\begin{equation}
	A_e = \tel{2E}\cdot e^{-\afrac{S}{K}}\cdot\int e^{\afrac{S}{K}}\cdot\intd{}J\,.
\end{equation}
There are only two cases where it is possible to directly compute the energy for state II, that is if either $S$ is constant or the material resists any change of volume. In those two cases the energy density with respect to the volume unit in the final state\Hcomment{\enquote{Final state} most probably refers to state II.} is \Herror{The number of the equation of was changed from (10c) to (10d).}
\begin{equation}
	A_e = \tel{2E}\cdot J\,.
\end{equation}
\end{subequations}
For constant $S$, the energy density with respect to the volume unit in the initial state is
\[
	A_a = \frac{e^{\afrac{S}{K}}}{2E}\cdot J\,,
\]
and for infinitely large K we find
\[
	A_a = \tel{2E}\cdot J\,,
\]
which was to be expected.\\
All these propositions hold only if the principal axes of deformation do not rotate, although in that case they hold for cylindrical coordinates and polar coordinates as well. We have not yet discussed the boundaries for the integrals of work. We must think of the $S_i$ and $e_i$ as functions in only one parameter, which we can interpret as the intensity of loading. Then the integral is to be taken from $0$ to a given value for this parameter. A different sequence of loading will generally result in a different amount of work\Hcomment{Hencky realizes that his law of elasticity is not hyperelastic. In today's notation his law can be written as\[\sigma = 2\,G\,\log V \;+\; \Lambda\,\tr[\log V]\cdot\id\,,\] where $G$ and $\Lambda$ are the two Lam\'e parameters and $\sigma$ denotes the Cauchy stress.}.

\subsection[The experiment of tension and compression for a cylindrical rod of an ideally elastic material]{The experiment of tension and compression\Hgerman{Zug- und Druckversuch} for a cylindrical rod of an ideally elastic material}
In \enquote{Elastizität und Festigkeit} \cite{von1905elastizitaet}, Bach describes the results of his experiments with rubber. His conclusion is that the elasticity modulus of rubber increases\Hcomment{Note that Bach uses the term \emph{modulus} to describe the \emph{inverse} of the quantities known as modulus today.} with extension and decreases with reduction, where stresses and extensions are defined as above. Since we assume fixed elastic constants, this can only be seen as pseudo-changes\Hgerman{scheinbare Änderungen}. To explain this phenomenon we will apply our equations to a very long cylindrical rod to which tension or compression is applied.\\
We choose a cylindrical coordinate system for this uniaxial state of stress. All functions are given with respect to the final state.\\
We denote by $r=\rbar-u$ the distance of a point to the axis and by $z=\zbar-w$ the original position of a point in direction of the axis.\\
Then
\begin{align*}
	\dr &= \drbar\cdot\left(1-\pdd{u}{\rbar}\right),\\
	\dz &= \dzbar\cdot\left(1-\pdd{w}{\zbar}\right),\\
	r\cdot\dphi &= \rbar\,\dphi\cdot\left(1-\frac{u}{\rbar}\right).
\end{align*}
Since for an infinitely long rod any cross section can be considered as positioned in the middle, $\pdd{w}{\zbar}$ must be independent of $\rbar$ for symmetry reasons and will be denoted by $\lambda$. For $u=-\rbar\cdot x$, where $x$ is a constant unknown for now, we find
\[
	\pdd{u}{\xbar} = -x \qquad \text{and} \qquad \frac{u}{\rbar} = -x\,.
\]
We will now identify the indices $1,2,3$ in equation \eqref{eq:Hencky4a} with $r$, $\varphi$ and $z$. We find $S_r = S_\varphi$ and $e_r=e_\varphi=-x,\ e_z=\lambda$. On the boundary, and hence on the whole cylinder, the condition $S_r=S_\varphi=0$ must hold, which implies
\[
	(1+x)^{k+\afrac23}\cdot(1+x)^{k-\afrac13}\cdot(1-\lambda)^{k-\afrac13}
=1\,,\]
thus after some transformations the ratio $x$ computes to
\[
	x=-1+(1-\lambda)^{-\frac{k-\afrac13}{2k+\afrac13}}\,.
\]
According to equation \eqref{eq:Hencky2a},
\[
	\tel{m} = \frac{k-\afrac13}{2k+\afrac13}
\]
and thus expanding $x$ into a series yields
\[
	x = \tel{m}\cdot\lambda + \tel{m}\cdot\left(\tel{m}+1\right)\cdot\frac{\lambda^2}{2}+\dots\,,
\]
from which we conclude that $m$ loses its meaning as the contraction coefficient once $\lambda$ is no longer very small. The stress $S_z$ computes to $S_z=2\,G\cdot\ln\frac{1+x}{1+\lambda}$. We insert the value of $x$, and with\Herror{$2\,G(1+\afrac1m)$ was changed to $E=2\,G(1+\afrac1m)$} $E=2\,G\,(1+\afrac1m)$ we obtain
\begin{subequations}
\begin{equation}
	S_z = -E\cdot\ln(1-\lambda)
\end{equation}
or, inside the convergence radius $1>\lambda>-1$,
\begin{equation}
	S_z = E\cdot\left(\lambda+\frac{\lambda^2}{2}+\frac{\lambda^3}{3}+\dots\right)\,.
\end{equation}
\end{subequations}
If we extend some piece to twice its initial length, i.e. $\lambda=\afrac12$, the elastic modulus determined through common means increases to around $1{.}4$ times its original value. For negative values of $\lambda$ we obtain a good correspondence to experimental results as well.\\
We can compute the elastic energy using equation \eqref{eq:Hencky10b}. With
\[
	3\cdot K = \frac{E\,m}{m-2}
\]
we find
\[
	A_a = \tel{E}\cdot\int e^{\afrac{S_z}{3K}}\cdot S_z\cdot \intd{}S_z\,.
\]
By integrating from $0$ to $S_z$ we obtain
\begin{equation}
	A_a = E\cdot\frac{m^2}{(m-2)^2}\cdot\left\{1-e^{\afrac{S_z}{3K}}\cdot\left(\frac{S_z}{3K}-1\right)\right\}\,.
\end{equation}
For very small values of $\afrac{S_z\,}{3K}$ we can expand into series and obtain the familiar value
\[
	A_a = \tel{2E}\cdot S_z^2\,.
\]
\subsection{The uniform expansion of an elastic membrane}
As a counterpart to the case considered above we will now assume $S_r=S_\varphi=\text{const.}$ and $S_z=0$. An example of such a state of stress is a thin plate being extended uniformly in its plane. We will give all functions with respect to the final state in this case as well and define
\[
	u=x\cdot\rbar \qquad \text{and} \qquad w=-\lambda\cdot\zbar\,,
\]
hence
\[
	\pdd{w}{\zbar} = -\lambda\,; \qquad \pdd{u}{\rbar}=x\,.
\]
Then the condition $S_z=0$ yields
\[
	(1+\lambda)((1-x)^2(1+\lambda))^{k-\afrac13}=1\,,
\]
which implies
\[
	\lambda = -1+(1-x)^{-\frac{2k-\afrac23}{k+\afrac23}}\,.
\]
The exponent can be easily expressed in terms of $m$,  we find
\[
	\frac{2k-\afrac23}{k+\afrac13} = \frac{2}{m-1}\,.
\]
The stress $S_r$ computes to $S_r = S_\varphi = 2\,G\ln[(1-x)^{k+\afrac23}(1-x)^{k-\afrac13}(1+\lambda)^{k-\afrac13}]$ or, after insertion of the value for $\lambda$ obtained above and some transformations,
\begin{equation}
	S_r=S_\varphi=-\frac{E\cdot m}{m-1}\cdot\ln(1-x)\,.\label{eq:Hencky13}
\end{equation}
The energy density with respect to the undeformed volume element is
\begin{align*}
	A_a &= \tel{E}\cdot\int e^{\afrac{S}{K}}\cdot\left\{\left(S_r-\tel{m}\cdot S_\varphi\right)\intd{}S_r \,+\, \left(S_\varphi - \tel{m}\cdot S_r\right)\intd{}S_\varphi\right\}\\
	&= \frac2E\left(1-\tel{m}\right)\cdot\int e^{\afrac{2S_r}{3K}}\cdot S_r\:\intd{}S_r\,,
\end{align*}
and integrating from $0$ to $S_r$ yields
\begin{equation}
	A_a = \frac{E\cdot m(m-1)}{2(m-2)^2}\cdot\left\{1+e^{\afrac{2S_r}{3K}}\cdot\left(\frac{2S_r}{3K}-1\right)\right\}\,.
\end{equation}
We can use the above formulae for a small application. If we inflate a thin-walled hollow sphere made from rubber, we can ignore the stress in the radial direction for a first order approximation. Then we can apply the above considerations, and we only have to find a relation between the excess inner pressure and the stresses $S_l$ in the elastic membrane.
If $h_x$ is the thickness of the membrane, $R_x$ is the radius of the inflated hollow sphere, and $h$ and $R$ are their respective initial values, then these quantities can be related with the stretches by
\[
	x=\frac{R_x-R}{R_x}\,, \qquad -\lambda=\frac{h_x-h}{h_x}
\]
and
\[
	h_x = h\cdot\left\{\frac{R}{R_x}\right\}^{\frac{2}{m-1}}\,,
\]
and equation \eqref{eq:Hencky13} yields
\[
	S_l = \frac{E\cdot m}{m-1}\cdot\ln\left\{\frac{R_x}{R}\right\}\,.
\]
Then we can obtain the relation between the pressure and the radius of the sphere from the equilibrium condition for a circular piece of the membrane in radial direction:
\[
	p = 2E\cdot\frac{h}{R}\cdot\frac{m}{m-1}\cdot\frac{\ln(\afrac{R_x}{R})}{(\afrac{R_x}{R})^{\frac{m+1}{m-1}}}\:.
\]
A range of simple special cases could be solved in this way. It would be particularly interesting to examine vibrations of finite amplitude, but at this point we will refrain from going into this in more detail. However, a comment on our attitude towards empirical observations seems appropriate. One could think that establishing a law of elasticity would be a matter of empirical research. This, however, is a misconception since there is no physical body which does not, for sufficiently large strains, exhibit plastic behaviour or change its original connectedness through  countless microscopic fissures. While it is certainly a matter of empirical observation to determine how actual materials compare to the ideally elastic body, the law itself acts as a measuring instrument which is extended into the realm of the intellect, making it possible for the experimental researcher to make systematic observations.\\
In this sense, whenever the strength of materials is considered beyond the scope of Hooke's law, experimental researchers require the assistance of the theorist as much as they require the help of craftsmen building their instruments in order to keep track and remain in control of their acquired observations\Hgerman{gewonnene Erfahrungen}.

\subsection*{Summary}
A law of elasticity valid for finite deformations is deduced. In doing so we show that the effect of an additional loading is independent of the already existing stress if the principal axes of the deformation remain aligned with the original axes of the coordinate system. We obtain a simple expression for the elastic energy as well as a simple theory for the tensile test of rubberlike materials, where the effect of permanent deformations and small fissures is ignored.

\newpage
\setcounter{equation}{0}
\setcounter{footnote}{0}
\section[Heinrich Hencky: What circumstances effect the hardening for ductile deformations of isotropic solid bodies?]{Heinrich Hencky: What circumstances effect the hardening for ductile deformations of isotropic solid bodies?\Hgerman{Welche Umstände bedingen die Verfestigung bei der bildsamen Verformung von festen isotropen Körpern? \cite{Hencky1929}}}

\dictum{We think of the material points of an elasto-plastic body as attached to two interpenetrating grids, one of which has large shear modulus and bulk modulus but exhibits stress states bound to the limit of plasticity, while the other one stays ideally elastic with a very small shear modulus and vanishing bulk modulus. By the example of the uniaxial stress state we show the form of the stress-tension-diagram for our model. Furthermore we show how one must distinguish between purely elastic, potentially elastic and lost energy in the elasto-plastic case.}

\subsection*{Introduction}
The research on the strength of materials suggests that the so-called hardening in plastic deformations is connected to the crystalline structure of quasi-isotropic materials. At least the hardening is an experimental fact, thus it is probably useful to consider the stress state of such a hardened material under the assumption of total isotropy, undeterred by the fact that such an ideally isotropic body does not exist.
\subsection{The law of elasticity for an ideally elastic body}
Before we are able to talk about deviations from ideally elastic behaviour, we must first establish the law of elasticity for finite deformations. While the assumption may be warranted that finite deformations are irrelevant because plastic deformation already takes place for very small displacements, we must consider finite deformations to properly justify our approach.\\
There are two conditions we impose on our ideally elastic material. First we require that the applied work is converted fully into elastic energy, which is released without loss after the loading is removed.\\
This first condition is fulfilled if there exists an elastic potential, but it allows for a number of different laws of elasticity. Thus we impose a second condition: if we apply a second loading to an already deformed body, it must not be possible to obtain the first loading from the deformations resulting from the additional loading.\\
This second condition is satisfiable only if the volume is constant, and it can only be satisfied in a single way.\\
It is this uniqueness of our law of elasticity which ranks our ideally elastic material as a mathematical ideal which is independent of empirical observation. Of course, the usefulness of such mathematical ideals is a matter of observation. The law of elasticity determines the relation between a stress tensor and the eigenvalues of an affine transformation matrix. Since the quantities used in the elementary theory of elasticity for the definition of strain are not symmetric with respect to the reference state and the final state\Hcomment{The strain measures do not satisfy the tension-compression symmetry $e(\lambda\inv) = -e(\lambda)$.}, they are not suitable for our purpose. Therefore, as a measure of the elastic deformation, we introduce the
\[
	\ln\left\{\frac{\text{final length}}{\text{initial length}}\right\}
\]
which agrees with the familiar measure of strain for very small deformations. For now we will consider only cases in which the principal axes of the deformation do not change. Our fixed coordinate system refers to the initial state. The side lengths of the volume element are
\begin{align*}
	&\text{$\dx_1,\,\dx_2,\,\dx_3 \quad $ in state I,}\\
	&\text{$\dxbar_1,\,\dxbar_2,\,\dx_3 \quad $ in state II,}\\
	&\text{$\dxbarbar_1,\,\dxbarbar_2,\,\dxbarbar_3 \quad $ in state III.}
\end{align*}
We introduce the displacements $u_i$ ($i=1,2,3$) and define
\begin{subequations}
\begin{align}
	\dxbari &= \dxi\left(1+\pdd{u_i}{x_i}\right) = \dxi(1+e_i) = \dxi\cdot e^{\eps_i}\,,\\
	\dxbarbari &= \dxi\left(1+\pdd{u_i}{x_i}+\pdd{\delta u_i}{x_i}\right) = \dxi(1+e_i+\delta e_i) = \dxi\cdot e^{\eps+\delta \eps_i}\,.
\end{align}
We obtain $\delta\eps_i = \pdd{\delta u_i}{x_i}$ for the conversion of the differentials as well as
\begin{equation}
	\dxbarbari = \dxbari\left(1+\frac{\delta e_i}{1+e_i}\right) = \dxbari\,(1+\delta\eps_i)
\end{equation}
for the transformation from II to III, which we assume to be infinitesimal. The change of volume is given by\Hcomment{Here, $\Delta=\det F$.}
\begin{equation}
	\frac{\dVbar}{\dV} = \Delta = (1+e_1)(1+e_2)(1+e_3) = e^{\eps_1+\eps_2+\eps_3}\,.
\end{equation}
\end{subequations}
The quantities $\eps$ are the logarithms of the affine ratios. The true physical stresses, with respect to the final state, are denoted by $S_i$. The mean stress is
\begin{subequations}
\begin{equation}
	S=\frac13\,\{S_1+S_2+S_3\}\,.
\end{equation}
Without yet knowing the law of elasticity, we are now able to give the differential of the performed work:
\[
	\delta A\cdot\dx_1\dx_2\dx_3 = S_1\dxbar_2\dxbar_3\cdot\pdd{\delta u_1}{x_1}\cdot\dx_1 + S_2\dxbar_3\dxbar_1\cdot\pdd{\delta u_2}{x_2}\cdot\dx_2 + S_3\dxbar_1\dxbar_2\cdot\pdd{\delta u_3}{x_3}\cdot\dx_3\,.
\]
Here the work is given with respect to the undeformed state; some easy computations yield
\begin{align}
	\delta A &= \Delta\left\{\frac{S_1\cdot\delta e_1}{1+e_1}+\frac{S_2\cdot\delta e_2}{1+e_2}+\frac{S_3\cdot\delta e_3}{1+e_3}\right\}\nnl
	&= \Delta\,\,\{S_1\,\delta\eps_1 + S_2\,\delta\eps_2 + S_3\,\delta\eps_3\}\,.\label{eq:29Hencky2b}
\end{align}
The occurrence of the change of volume in this expression is due to the fact that the stress is not actually a tensor, but a tensor density. Thus the principal stresses are scalar densities. Since the eigenvalues of the transformation matrix are true invariants, the law of elasticity cannot be expressed in terms of the quantities $S_i$. Instead we must provide an absolute invariant\Hcomment{Since $S_i$ denote the principal Cauchy stresses, the quantities $T_i = \det(F)\cdot S_i$ are the \emph{principal Kirchhoff stresses}.}
\begin{equation}
	T_i = \Delta\cdot S_i
\end{equation}
\end{subequations}
for which we can write down the law of elasticity in a linear form. With $\eps=\frac13(\eps_1+\eps_2+\eps_3)$ we find\Hcomment{At this point, Hencky introduces a law of elasticity which differs from the law he proposed in his 1928 article: instead of the Cauchy stress, he relates the \emph{Kirchhoff stress} $\tau$ to the stretch tensor $V$ and obtains\[\tau = 2\,G\,\log V \;+\; \Lambda\,\tr[\log V]\cdot\id\,.\]}
\begin{subequations}
\begin{equation}
	T_i=2\,G\cdot\{\eps_i+(3k-1)\eps\}
\end{equation}
and thus
\begin{equation}
	T=2\,G\,k\cdot3\eps=3K\cdot\eps\,.
\end{equation}
Then the law of superposition is linear as well:
\begin{equation}
\label{eq:29Hencky3c}
	\delta T_i = 2\,G\cdot\{\delta\eps_i + (3k-1)\delta\eps\}\,.
\end{equation}
\end{subequations}
If we insert the values for $\delta\eps_i$ from equations \eqref{eq:29Hencky3c} into equations \eqref{eq:29Hencky2b}, we obtain a total differential $\delta A$ and, by integrating, the stored elastic energy:\Hcomment{The energy function obtained from Hencky's new law of elasticity is known today as the \emph{quadratic Hencky strain energy}:\[W(V) = G\,\norm{\dev_3 \log V}^2 + \frac{K}{2}\,[\tr(\log V)]^2\,.\]}
\begin{subequations}
\begin{align}
	2A &= \tel{2\,G}\cdot\{(T_1-T)^2+(T_2-T)^2+(T_3-T)^2\}+\tel{K}\cdot T^2\,,\\
	2A &= \frac{\Delta^2}{2\,G}\{(S_1-S)^2+(S_2-S)^2+(S_3-S)^2\}+\frac{\Delta^2}{K}\cdot S^2
\intertext{or finally, with respect to the deformation components,}
	2A &= 2\,G\cdot\{(\eps_1-\eps)^2+(\eps_2-\eps)^2+(\eps_3-\eps)^2\}+9\,K\eps^2\,.
\end{align}
\end{subequations}
If the principal axes of the deformation do not remain fixed, we can still use the same expression to obtain formulae for the general spatial problem.

\newpage
\setcounter{equation}{0}
\setcounter{footnote}{0}

\makeatletter
\tagsleft@true
\let\veqno\@@leqno
\makeatother 

\section[Heinrich Hencky: The law of superposition for a finitely deformed elastic continuum capable of relaxation and its significance for an exact derivative of the equations for the viscous fluid in the Eulerian form]%
{Heinrich Hencky: The law of superposition for a finitely deformed elastic continuum capable of relaxation and its significance for an exact derivative of the equations for the viscous fluid in the Eulerian form\Hgerman{Das Superpositionsgesetz eines endlich deformierten relaxationsfähigen elastischen Kontinuums und seine Bedeutung für eine exakte Ableitung der Gleichungen für die zähe Flüssigkeit in der Eulerschen Form \cite{Hencky1929super}}}

\dictum{A generally nonlinear law of superposition is deduced from the law of elasticity for an ideally elastic material. Using the approach of a relaxation in the Maxwellian sense we develop the basic equations for the motion of viscous fluids from the finite deformation of an initially elastic body. We show that that the Navier-Stokes equations of hydrodynamics emerge as a special case for an infinitely large shear modulus and an infinitely small relaxation time, where the product of these two quantities remains finite, but that they must fail for inhomogeneous flow states of real fluids if large velocities are attained.}

\subsection*{Introduction}
Until a few decades ago the viscous liquid was interpreted, especially from the physical point of view, as an elastic body, where the deviatoric part of the strain energy is subject to a quick absorption by the thermal motion of small particles \cite{natanson1901a, natanson1901b, reiger1919} (c.f. \cite[p. 102]{auerbach1931handbuch} including the references given there; strangely enough, the Maxwellian theory is not even mentioned in \cite{grammel1927mechanik}).

This derivation of the hydrodynamic equations has been developed in a number of articles, R. Reiger in particular has repeatedly advocated the Maxwellian interpretation. If, in the following, we have to criticize these previous works for their accomplishment with respect to the theory of elasticity, it is only to free the valuable ideas laid out in these works from their obscuring ambiguities\Hgerman{hindernde Unklarheiten}. The existence of such ambiguities is shown by the fact that the theory of relaxation, as opposed to the concept of Newtonian friction, has had no influence on the development of hydrodynamics.

Another point, of course, is that the significance of the model of relaxation has been overestimated at first. Experiments show that at least polycrystalline materials \emph{below} their limit of elasticity are not subjected to relaxation, and it is pointless\Hgerman{zwecklose Spielerei} to attribute to them an infinite relaxation time. Apart from that, it seems premature to ignore such a prolific and original concept of friction in viscous liquids. The Maxwellian conception yields an intuitively reasonable concept\Hgerman{Gedankenmodell\label{footnote:Gedankenmodell}} and establishes a mechanical transition from the solid to the liquid continuum, the epistemological justification of which is beyond doubt, and its deductions remain to be properly examined through experiments.

Here, our aim is to show that even a purely theoretical analysis of Maxwell's concept\Href{footnote:Gedankenmodell} of the viscous fluid is able to prove the necessity of revising the currently prevalent interpretation of friction for large flow rates.

\subsection{The law of elasticity for an ideally elastic material}
In one of the articles in the A.Föppl commemorative volume \cite{prandtl1924}, L. Prandtl considers a number of constructions which he divides into elastically determinate and elastically indeterminate ones. By Prandtl's definition, elastically determinate constructs are those for which the changes induced by the application of an additional loading are independent of already occurring stresses. If, however, the deformation depends not only on the additional loading but also on internal stresses of the system in the previous state, Prandtl calls the construct elastically determinate.

The underlying concept of Prandtl's distinction is that of an algebraic group, since for an elastically determinate continuum the group of changes of shape must be isomorphic to the group of changes of stress.

The law of elasticity requires a relation between the eigenvalues of a stress tensor and those of a transformation matrix which is as simple as possible. However, the solely considered pure deformations, i.e. those transformations which can be expressed as three stretches along mutually orthogonal directions, do not form a transformation group; the pure deformations form a group only if the principal axes of deformation do not rotate. The increments of stress, which in this special case are equal to the transformations, must form an isomorphic group and are therefore independent of earlier increments, thus we find an elastically determinate law of superposition in Prandtl's sense. If, however, the principal axes of deformation rotate, the corresponding pure deformations lose their group properties\Hgerman{Gruppencharakter} and only an elastically indeterminate law of superposition is possible.

Now, in setting up a law of elasticity for an ideally elastic material, we require elastic determinacy to the greatest extent for epistemological reasons (c.f. \cite{dingler1928experiment} for a more general explanation). Just like the rigid body, the ideally elastic material is not a real material but an instrument of measurement and comparison, therefore it must be left to the experiment to discover elastic indeterminacies which can be avoided in the theory.

There are two reasons why elastic indeterminacy cannot be avoided completely. The first one, the lack of group properties for pure deformations in the general case, has already been discussed. The second reason is that the stress tensor is not a true tensor of weight 0 but a tensor density (this was first suggested by L. Brillouin, c.f. \cite{brillouin1925}).

In the elementary theory of elasticity, these issues have been neglected so far under the reasoning that only very small deformations occur; it has been overlooked that certain important geometrical relations, which could not have been found through observation\Hgerman{Erfahrung} alone, are blurred and obfuscated by the transition to the infinitesimal case, when they should actually constitute the foundation of the mathematical theory of elasticity.

Since tensor densities obtain the transformation properties of tensors only through multiplication with the determinant of the transformation, we must distinguish the true physical stress $S_i$ from the reduced stress $S_i' = ~$volume$\,\times S_i$ which a law of elasticity must relate to a function of the principal strains.

As we have shown in a number of earlier works \cite{Hencky1928, Hencky1929, biezeno1929}, the conditions of elastic determination mentioned above are satisfied if we measure the strain by $\ln\left\{\frac{\text{final length}}{\text{initial length}}\right\}$.

Denoting the principal strains by $\eps_i$ we can formulate the elastic energy for finite deformations even for changing directions of the principal axes of deformation:
\begin{equation}
	2A = 2G\cdot \{(\eps_1-\eps)^2 + (\eps_2-\eps)^2 + (\eps_3-\eps)^2\} + 9K\cdot \eps^2\,.\label{eq:a}
\end{equation}
The law of elasticity corresponding to this approach, with
\[
	\eps = \tel3 (\eps_1 + \eps_2 + \eps_3)
\]
and
\[
	S_i' = e^{3\eps} \cdot S_i\,,
\]
is
\begin{subequations}
\label{eqs:superHencky2}
\begin{align}
	S_i' - S' &= 2G\cdot \{\eps_i - \eps\}\label{eq:superHencky2a}\\
	S' &= 3K\cdot\eps\,.\label{eq:superHencky2b}
\end{align}
\end{subequations}
For very elastic bodies, e.g. rubber, this law has indeed been confirmed through experiments.

As a distinction to other fields of physics we will henceforth assume isothermal changes of state, possibly generated heat will be thought of as automatically dissipated.

The law we assumed is only the simplest of all possible laws, therefore one should not think that the far reaching conclusions we will draw later on are founded in the special form of these equations. In this regard we could also have based our work on the law used in elementary elasticity which, for very small deformations, can be obtained from equations \eqref{eqs:superHencky2} as well.

\subsection{The law of superposition for the ideally elastic continuum}
The transition towards an arbitrary system of coordinate axes not coinciding with the principal axes involves some difficulties which do not occur in the case of infinitesimal deformations. Under the assumption of an arbitrarily oriented system we must therefore look for an analogy of equations \eqref{eq:superHencky2a} and \eqref{eq:superHencky2b}.

To participate in the advantages of tensor notation for the necessary computations, it is by no means necessary to employ a skew coordinate system. It suffices to number the coordinate directions and, unless otherwise indicated, sum over indices occurring twice without writing an explicit summation sign.

We distinguish three different states:
\begin{align*}
	&\text{state \;\,0 with the coordinates $\xring_i = x_i-u_i$\,;}\\
	&\text{state \;\,I\, with the coordinates $x_i$\,;}\\
	&\text{state II, initially with the coordinates}\\
	&\qquad\qquad \xbar_i = x_i + \delta u_i = x_i + v_i \, \delta t\,,
\end{align*}
where
\[
	v_i = \frac{\delta u_i}{\delta t}\,.
\]
Thus the system is based on state I, while states 0 and II are derived from it.

We must now aim to find a quantity which completely describes the given state of deformation. Through total differentiation we obtain:
\begin{subequations}
\begin{align}
	d\xring_i = \pdd{\xring_i}{x_1}\cdot dx_1 + \pdd{\xring_i}{x_2}\cdot dx_2 + \pdd{\xring_i}{x_3}\cdot dx_3 &= \pdd{\xring_i}{x_k}\cdot dx_k\,,
\intertext{as well as the identity:}
	dx_i &= \pdd{x_i}{x_k}\cdot dx_k\,. \nonumber
\end{align}
Using our summation convention we can now write:
\begin{align*}
	dx_i \cdot dx_i &= dx_i^2 = \pdd{x_i}{x_k} \cdot \pdd{x_i}{x_l} \cdot dx_k\,dx_l\,,\\
	d\xring_i \cdot d\xring_i &= d\xring_i^2 = \pdd{\xring_i}{x_k} \cdot \pdd{\xring_i}{x_l} \cdot dx_k\,dx_l\,,
\end{align*}
where we sum over nine terms on each right hand side, and subtracting these two quantities yields a measure of the pure deformation, eliminating the rotational movement of the volume element which is not relevant to us.

To define the tensor quantity $e_{kl}$, which is characteristic of the pure deformation\Hcomment{Note carefully that $e_i$ does not denote the same quantity as in the previous article: here, $e_i = \half(1-\frac1{\lambda_i^2})$ are the eigenvalues of the strain tensor $\half(\id-C\inv)$.}, we let
\begin{equation}
	dx_i^2 - d\xring_i^2 = 2\cdot e_{kl} \cdot dx_k\,dx_l
\end{equation}
\end{subequations}
and find:
\[
	2\cdot e_{kl} = \pdd{x_i}{x_k}\cdot\pdd{x_i}{x_l} - \pdd{\xring_i}{x_k}\cdot\pdd{\xring_i}{x_l}
\]
or, after some computation:
\begin{subequations}
\begin{equation}
	2\cdot e_{kl} = \pdd{u_l}{x_k} + \pdd{u_k}{x_l} - \pdd{u_i}{x_k}\cdot\pdd{u_i}{x_l}\,.
\end{equation}
The eigenvalues of this tensor can be obtained from the cubic equation\Hcomment{Here, the term $|A|$ denotes the determinant of a matrix $A$.}:
\begin{equation}
	\left|\begin{array}{lll}
		e_{11}-e_i\: & e_{12} & e_{13}\\
		e_{21} & e_{22}-e_i\: & e_{23}\\
		e_{31} & e_{32} & e_{33}-e_i
	\end{array}\right| = 0\,.
\end{equation}
There is a simple relation between these eigenvalues and the quantities $\eps_i$ we used in equations \eqref{eqs:superHencky2}. The conversion into the values $e_i$ can be easily accomplished by means of the equation defining the $e_{kl}$. We find:
\begin{equation}
	\left\{\begin{aligned}
		\eps_i &= \ln\left\{\frac{1}{\sqrt{1-2\,e_i}}\right\}\\
		2\,\eps_i &= -\ln(1-2\,e_i)\,.
	\end{aligned}\right.
\end{equation}
\end{subequations}
This form does not yet allow us to consider an arbitrary coordinate system, we first need to express the relation in the form of a series converging absolutely for all values of $e_i$. We obtain:
\begin{subequations}
\begin{equation}
	2\,e_i = 1-e^{-2\eps_i} = \frac1{1!}\cdot(2\eps_i) - \frac1{2!}\cdot(2\eps_i)^2 + \frac1{3!}\cdot(2\eps_i)^3 - \dots
\end{equation}
On the other hand we can use equations \eqref{eqs:superHencky2} to express the quantities $\eps_i$ in terms of $S_i'$, namely:
\begin{equation}
	2\eps_i = \frac1G \cdot (S_i'-S') + \frac2{3\,K}\cdot S'\,.
\end{equation}
By introducing the shortened notation:
\begin{equation}
	\left\{\begin{aligned}
		\quad\sigma_i' &= \frac1G \cdot (S_i'-S') \quad\text{and}\\
		\sigma' &= \frac2{3\,K}\cdot S'\,,
	\end{aligned}\right.
\end{equation}
\end{subequations}
we obtain the important equality:
\begin{equation}
	2\cdot e_i = \frac1{1!}\cdot(\sigma_i'+\sigma') - \frac1{2!}\cdot(\sigma_i'-\sigma')^2 + \frac1{3!}\cdot(\sigma_i'+\sigma')^3 - \dots\label{eq:superHencky6}
\end{equation}
which provides the basis for all further examinations. If we had chosen a different law of elasticity, this relation would have become even more complicated. Indeed, our choice was made such that the relation takes on the simplest form possible.

For the transition from powers, such as $e_i^1,e_i^2,e_i^3$, referring to the principal axes to an arbitrary orientation of coordinates, our notation of indices provides an easy method, which we will state without proof since it can be found in any textbook on algebra.

For example, if $e_i^3$ is given and $e_{m\,n}$ is the tensor form of $e_i$, then the general tensor form of $e_i^3$ can be represented as $e_{m\,i}\,e_{i\,k}\,e_{k\,n}$, where:
\begin{align*}
	e_{m\,i}\,e_{i\,k}\,e_{k\,n} &= e_{m\,1}\,e_{1\,k}\,e_{k\,n} + e_{m\,2}\,e_{2\,k}\,e_{k\,n} + e_{m\,3}\,e_{3\,k}\,e_{k\,n}\\
	&= e_{m\,1}\,e_{1\,1}\,e_{1\,n} + e_{m\,1}\,e_{1\,2}\,e_{2\,n} + e_{m\,1}\,e_{1\,3}\,e_{3\,n}\\
	&\,+ e_{m\,2}\,e_{2\,1}\,e_{1\,n} + e_{m\,2}\,e_{2\,2}\,e_{2\,n} + e_{m\,2}\,e_{2\,3}\,e_{3\,n}\\
	&\,+ e_{m\,3}\,e_{3\,1}\,e_{1\,n} + e_{m\,3}\,e_{3\,2}\,e_{2\,n} + e_{m\,3}\,e_{3\,3}\,e_{3\,n}\,.
\end{align*}
Using this simple pattern it is always possible to move on from the principal axes of a tensor to the general component representation.

We now pass on to state II, meaning we subject state I to an affine spatial transformation defined through the infinitesimal translation $\delta u_i$\,. We can decompose the transformation matrix into a symmetric part, the change of shape
\begin{subequations}
\begin{equation}
	f_{m\,n}\,\delta\,t = \half\cdot \left\{\pdd{\,\delta u_n}{x_m} + \pdd{\,\delta u_m}{x_n}\right\} = \half\cdot \left\{\pdd{v_n}{x_m} + \pdd{v_m}{x_n}\right\}\,\delta t
\end{equation}
and an antisymmetric part, the rotation:
\begin{equation}
	\omega_{m\,n}\,\delta\,t = \half\cdot \left\{\pdd{\,\delta u_n}{x_m} - \pdd{\,\delta u_m}{x_n}\right\} = \half\cdot \left\{\pdd{v_n}{x_m} - \pdd{v_m}{x_n}\right\}\,\delta t\,.
\end{equation}
\end{subequations}
Then $f_{m\,n}$ is the tensor of the deformation velocities and $\omega_{m\,n}$ is the tensor of angular velocities.

To find the corresponding changes to the tensor $e_{m\,n}$, we note that the rotation of the element must remain without effect on the law of superposition we are trying to find.

Thus the variational symbol refers to the change of shape of the material element we put our focus on.

However, we must not write:
\[
	e_i + \delta\,e_i = e_i + f_i\, \delta\,t
\]
or:
\[
	e_{m\,n} + \delta\,e_{m\,n} = e_{m\,n} + f_{m\,n}\cdot\delta\,t\,,
\]
since the infinitesimal spatial transformation slightly changes the tensor $e_{m\,n}$ as well. When computing this change we must consider the fact that we always compute the principal axes with respect to the final state and must therefore transform the tensor in the opposite direction.

According to the rules of tensor calculus each index is considered a representative of an ideal vector and is transformed as such. Then, if we introduce the identity tensor $g_{m\,n}$ with the matrix representation
\[
	\begin{matrix}
		1&0&0\\0&1&0\\0&0&1&\!\!\!\!\,,
	\end{matrix}
\]
we find
\[
	e_{m\,n} + \delta\,e_{m\,n} = (e_{i\,k} + f_{i\,k}\,\delta\,t)\,(g_{m\,i}-f_{m\,i}\,\delta\,t)\,(g_{n\,k}-f_{n\,k}\,\delta\,t)
\]
and after computation:
\begin{subequations}
\begin{equation}
	\delta\,e_{m\,n} = \delta\,t\,(f_{m\,n} - e_{n\,i}\,f_{i\,m} - e_{m\,k}\,f_{k\,n})\,.\label{eq:superHencky8a}
\end{equation}
From equations \eqref{eq:superHencky6} and \eqref{eq:superHencky8a} we finally obtain the desired law of superposition for a previous finite deformation of the material with which we master every problem in the theory of elasticity.
\begin{align}
	&\left\{\begin{aligned}
		2\cdot\{&f_{m\,n} - e_{n\,i}\,f_{i\,m} - e_{m\,k}\,f_{k\,n}\}\\
		&= \tel{1!}\cdot \left\{ \frac{\delta\,{\sigma_{m\,n}}'}{\delta\,t} + g_{m\,n}\cdot\frac{\delta\,\sigma'}{\delta\,t}\right\}\\
		&- \tel{2!}\cdot \left\{ \frac{\delta}{\delta\,t}({\sigma_{m\,i}}'\cdot{\sigma_{i\,n}}') + \frac{\delta}{\delta\,t}(2\cdot\sigma'\cdot{\sigma_{m\,n}}') + g_{m\,n}\,\frac{\delta}{\delta\,t}(\sigma')^2\right\}\\
		&+ \tel{3!}\cdot\{\dots\} \:-\: \dots
	\end{aligned}\right.
\intertext{The deformation quantities $e_{m\,n}$ occurring here are eliminated through the series:}
	&\left\{\begin{aligned}
		2\,e_{m\,n} = \tel{1!}\, ({\sigma_{m\,n}}' + g_{m\,n}\cdot\sigma') &- \tel{2!}\cdot({\sigma_{m\,i}}'\cdot{\sigma_{i\,n}}' + 2\cdot\sigma'\cdot{\sigma_{m\,n}}'\\
		&+ g_{m\,n}\cdot\sigma'^{\,2}) + \dots
	\end{aligned}\right.
\end{align}
\end{subequations}
After this elimination, the change of the stress tensor is given as a function of the tensor of deformation velocities and the tensor of velocity vectors. \emph{The displacements $u_i$, on the other hand, have disappeared from our formulae completely.}

The validity of these formulae is not limited with respect to the size of the ${\sigma_{m\,n}}'$\,.

If the ${\sigma_{m\,n}}'$ can be considered very small compared to unity, for example in the case of small elastic oscillations about a state of equilibrium, we can ignore higher powers of ${\sigma_{m\,n}}'$ and obtain the law of superposition in the following simpler form:
\begin{subequations}
\begin{align}
	\frac{\delta\,{\sigma_{m\,n}}'}{\delta\,t} &= \pdd{\,v_m}{\,x_n} + \pdd{\,v_n}{\,x_m} + \frac23\cdot g_{m\,n}\cdot\pdd{\,v_i}{\,x_i}\\
	\frac{\delta\,\sigma'}{\delta\,t} &= \frac23\cdot\pdd{\,v_i}{\,x_i}\,.
\end{align}
\end{subequations}
Now these equations are, without further thought, applied to finite movements and taken as a basis to deduce the equations of hydrodynamics for the theory of relaxation \cite{madelung1964mathematischen}.

The existence of such a linear law of superposition for arbitrarily large movements is only possible if the pure deformations form a transformation group. However, this is the case \emph{only} if the principal axes do not rotate at all.

We have yet to specify what is meant by variation of the stress state. The stresses are given as positional functions of the time $t$ for state I.

When forming the differential we must follow our axial trihedron\Hgerman{Achsendreikant} attached to the material particle and consider the rotation of the trihedron as well.

We therefore obtain:
\begin{equation}
	\frac{\delta\,{\sigma_{m\,n}}'}{\delta\,t} = \pdd{\,{\sigma_{m\,n}}'}{\,t} + v_i\cdot\pdd{\,{\sigma_{m\,n}}'}{x_i} + \,{\sigma_{n\,i}}'\cdot\omega_{i\,m} + \,{\sigma_{m\,k}}'\cdot\omega_{k\,m}\,.
\end{equation}

\newpage
\renewcommand{\refname}{References}

\end{document}